# Time continuity of weak-predictable random field solutions
Dr. Ejighikeme McSylvester **Omaba**

Department of Mathematics/Computer Science/ Statistics/Informatics
Faculty of Science
Federal University Ndufu-Alike Ikwo, Ebonyi State, Nigeria.
Email: mcsylvester.omaba@funai.edu.ng;mcsylvester_omaba@yahoo.co.uk

## ABSTRACT

The question of global existence or non-existence of solution to a given stochastic partial differential equation under some non-linear conditions always comes to mind. To show that our weak-predictable random field solutions do not have global existence for all time $t$, it requires that we first establish that the solutions exhibit continuity in time property. The results discuss the mean-square and mean continuity in time of a class of jump-discontinuous heat equations perturbed by compensated and non-compensated Poisson random noises respectively; and we showed that our mild solutions are mean and mean-square continuous in time for any time interval $[t_1, t_2], t_1, t_2 \geq 0$; better put, our solutions have continuous versions or modifications for any time interval.

## 1 INTRODUCTION

Consider the following stochastic heat equations driven by a compensated Poisson noise and non-compensated Poisson noise respectively,

$$[\frac{\partial u}{\partial t}(t,x) - \mathsf{L}u(t,x)]\mathrm{d}x\mathrm{d}t = \lambda \int_{\mathbf{R}} \sigma(u(t,x),h)\tilde{N}(\mathrm{d}h,\mathrm{d}x,\mathrm{d}t), \tag{1.1}$$

and

$$[\frac{\partial u}{\partial t}(t,x) - \mathsf{L}u(t,x)]\mathrm{d}x\mathrm{d}t = \lambda \int_{\mathbf{R}^d} \sigma(u(t,x),h)N(\mathrm{d}h,\mathrm{d}x,\mathrm{d}t), \tag{1.2}$$

with initial condition $u(0,x) = u_0(x)$. Here and throughout, $u_0 : \mathbf{R} \to \mathbf{R}_+$ is a non-random bounded and measurable function, and $\mathsf{L}$ is the $L^2$- generator of a real-valued Lévy process $X_t$ with the symbol $\Psi$ such that $E[e^{iX_t}] = e^{-t\Psi(\xi)}$, for every $\xi \in \mathbf{R}$, $t \geq 0$ and for all measurable function $f : \mathbf{R} \to \mathbf{R}_+$ satisfying

$$\mathrm{E}^x[f(X_t)] = \int_{\mathbf{R}} p(t,x,y)f(y)\mathrm{d}y = \int_{\mathbf{R}} f(y)\mathrm{P}^x(X_t \in \mathrm{d}y).$$

The above equations studied here are discontinuous analogues of the equations introduced in (Foondun and Khoshnevisan, 2009; Foondun and Khoshnevisan, 2013) and also more general than those considered in (Bie, 1998). Our mild solutions to these equations are weak-predictable random field solutions to the class of stochastic heat equations with jump-type process (Poisson random noise measure). Some precise conditions for existence and uniqueness of the solutions were given below and we had consequently established that the solutions grow in time at most a precise exponential rate at some time interval see (Omaba, 2014); and if the solutions satisfy some non-linear conditions then they cease to exist at some finite time $t$, see (Omaba et al., 2017). Therefore we study the mean and mean-square continuity in time of the weak-predictable random field solutions with non-compensated and compensated Poisson random measure respectively. This paper was motivated by works initiated by Khoshnevisan and his coauthors, see (Foondun and Khoshnevisan, 2009; Foondun and Khoshnevisan, 2013) and their references; where they systematically studied the large time behaviour of some class of stochastic equations. Lévy noise $\tilde{N}(dt,dx,dh)$ or $N(dt,dx,dh)$ has better modelling characteristics (Carr et al., 2002; Carr et al., 2003); the Lévy-type perturbations produce a better modelling result and performance of those natural occurrences and phenomena of some real world modelling, capturing some large moves and unpredictable events unlike Brownian motion perturbation that has many imperfections. The Lévy noise has a very rich and vast applications in Finance, Economics, Physics, etc.

## 2   UNDERLYING CONCEPTS, DEFINITIONS AND THEOREMS

**Definition 2.1** *We say that a process* $\{u(t,x)\}_{x\in\mathbf{R},t>0}$ *is a mild solution of (1.1) if a.s, the following is satisfied*

$$u(t,x) = \int_{\mathbf{R}} p(t,x,y)u_0(y)dy + \lambda \int_0^t \int_{\mathbf{R}} \int_{\mathbf{R}} p(t-s,x,y)\sigma(u(s,y),h)\tilde{N}(dh,dy,ds), \qquad (2.1)$$

where $p(t,.,.)$ is the heat kernel, see (Albeverio et al., 1998) and (Nicolas, 2000) and their references. If in addition to the above, $\{u(t,x)\}_{x\in\mathbf{R},t>0}$ satisfies the following condition

$$\sup_{0\leq t\leq T}\sup_{x\in\mathbf{R}} E|u(t,x)|^2 < \infty,$$

for all $T > 0$, then we say that $\{u(t,x)\}_{x\in\mathbf{R},t>0}$ is a random field solution to (1.1).

Now define

$$\Upsilon(\beta) := \frac{1}{2\pi}\int_{\mathbf{R}} \frac{d\xi}{\beta + 2\mathrm{Re}\Psi(\xi)} \qquad \text{for all } \beta > 0,$$

where $\Psi$ is the characteristic exponent for the Lévy process, otherwise known as the Lévy exponent. A result of (Dalang, 2004) shows that equation (1.1) has a unique solution with the requirement that $\Upsilon(\beta) < \infty$ for all $\beta > 0$ which forces $d = 1$ and coincides to a similar situation to that in (Foondun and Khoshnevisan, 2009) since $\tilde{N}$ is a martingale-valued Poisson measure and the integral in the definition of mild solution defined in terms of (Walsh, 1986). Fix some $x_0 \in \mathbf{R}$ and define the *upper $p$ th-moment Liapunov exponent* $\bar{\gamma}(p)$ of $u$ [at $x_0$] as

$$\bar{\gamma}(p) := \limsup_{t\to\infty} \frac{1}{t} \ln E\left[|u(t,x_0)|^p\right] \quad \text{for all } p \in (0,\infty), \qquad (2.2)$$

We assume the following condition on $\sigma$; which says essentially that $\sigma$ is globally Lipschitz in the first variable and bounded by another function in the second variable:

**Condition 2.2** *There exist a positive function $J$ and a finite positive constant, $\text{Lip}_\sigma$ such that for all $x, y, h \in \mathbf{R}$, we have*

$$|\sigma(0,h)| \leq J(h) \quad \text{and} \quad |\sigma(x,h) - \sigma(y,h)| \leq J(h)\text{Lip}_\sigma |x-y|.$$

The function $J$ is assumed to satisfy the following integrability condition:

$$\int_{\mathbf{R}} J(h)^2 \nu(dh) \leq K,$$

where $K$ is some finite positive constant.

**Remark 2.3** *Unlike the compensated noise term $\tilde{N}(dt, dx, dh)$, the non-compensated noise $N(dt, dx, dh)$ is not a martingale-valued Poisson random measure. The existence and uniqueness of the solution to (1.1) does not depend on the integrability condition (2.2) and the first moment of the solution exists; hence the existence and uniqueness for all $d \geq 1$.*

**Definition 2.4 (of random field solution)**
We seek a mild solution to equation (1.2) of the form.

$$u(t,x) = \int_{\mathbf{R}^d} p(t,x,y) u_0(y) dy + \lambda \int_0^t \int_{\mathbf{R}^d} \int_{\mathbf{R}^d} p(t-s,x,y) \sigma(u(s,y),h) N(dh, dy, ds), \tag{2.3}$$

with $p(t,.,.)$ the heat kernel. We impose the following integrability condition on the solution.

$$\sup_{t>0} \sup_{x \in \mathbf{R}^d} \mathrm{E}|u(t,x)| < \infty.$$

Let us define a Poisson random measure $N = \sum_{i \geq 1} \delta_{(T_i, X_i, Z_i)}$ on $\mathbf{R}_+ \times \mathbf{R}^d \times \mathbf{R}^d$ defined on a probability space $(\Omega, \mathsf{F}, P)$ with intensity measure $dt dx \nu(dh)$ where $\nu$ is a Lévy measure on $\mathbf{R}^d$; that is, it satisfies the following

$$\int_{\mathbf{R}^d} (1 \wedge h^2) \nu(dh) < \infty.$$

Here, the following assumptions were made on $\sigma$

**Condition 2.5** *There exist a positive function $J$ and a finite positive constant, $\text{Lip}_\sigma$ such that for all $x, y, h \in \mathbf{R}^d$, we have*

$$|\sigma(0,h)| \leq J(h) \quad \text{and} \quad |\sigma(x,h) - \sigma(y,h)| \leq J(h)\text{Lip}_\sigma |x-y|.$$

The function $J$ is assumed to satisfy the following integrability condition.

$$\int_{\mathbf{R}^d} J(h) \nu(dh) \leq K,$$

where $K$ is some finite positive constant.

**Definition 2.6** *(Solution norms.)* We define the norm on the mild solution (2.1)

$$\|u\|_{2,\beta} := \left\{ \sup_{t>0} \sup_{x \in \mathbf{R}} e^{-\beta t} \mathrm{E}[|u(t,x)|^2] \right\}^{1/2}$$

Likewise, we give the first moment norm of the solution (2.3) as follows:

$$\|u\|_{1,\beta} = \sup_{t \geq 0} \sup_{x \in \mathbf{R}^d} e^{-\beta t} \mathrm{E}|u(t,x)|.$$

We will apply Kolmogorov's continuity theorem to establish the time continuity of our weak-predictable random field solution and consider the increment $E[|u(t+h,x)-u(t,x)|^p]$ for $h \in (0,1)$ and $p=1,2$. We begin with the following definitions. There are different concepts or characterisations of continuity that appear in literature:

**Definition 2.7** *A function $f : \mathbf{R}^+ \to \mathbf{R}$ is called locally $\gamma$-Hölder continuous at $t \geq 0$ if there exist $\varepsilon, \gamma > 0$ and $C = C(t) > 0$ such that*
$$|f(t)-f(s)| \leq C|t-s|^\gamma,$$
for every $s \geq 0$ such that $|t-s| < \varepsilon$. The number $\gamma$ is called Hölder exponent and $C = C(t)$ is the Hölder constant. The function is Lipschitz continuous if $\gamma = 1$.

**Proposition 2.8** *(Kolmogorov-$\hat{C}$entsov Theorem). Let $(\Omega, \mathbf{F}, \mathbf{P})$ be a probability space and $(X_t)_{t \geq 0}$ a real-valued random process on it. Suppose there exist $\alpha \geq 1, \beta, C > 0$ such that*
$$E_\mathbf{P}[|X_t - X_s|^\alpha] \leq C|t-s|^{\beta+1}.$$
Then there exists a version $(Y_t)_{0 \leq t \leq T}$ of $(X_t)_{0 \leq t \leq T}$ for every $T \geq 0$ with $\gamma$-Hölder continuous paths for any $0 \leq \gamma \leq \beta/\alpha$.

In this paper, we will discuss the mean and the mean-square continuity in time of our random field solutions under some conditions on $\sigma$. The following definitions were informed by the above Proposition 1.8.

**Definition 2.9** *The mild solution (2.1) is said to be mean-square continuous in time if $t_1 < t_2$ such that*
$$\lim_{\delta \downarrow 0} \sup_{|t_1-t_2|<\delta} E|u(t_2,x)-u(t_1,x)|^2 = 0, \text{ for a fixed } x \in \mathbf{R}$$
or
$$\lim_{t_1 \uparrow t_2} E|u(t_2,x)-u(t_1,x)|^2 = 0, \text{ for a fixed } x \in \mathbf{R}.$$
Similarly, the mild solution (2.3) is said to be mean continuous in time if $t_1 < t_2$ and
$$\lim_{\delta \downarrow 0} \sup_{|t_1-t_2|<\delta} E|u(t_2,x)-u(t_1,x)| = 0, \text{ for a fixed } x \in \mathbf{R}^d$$
or
$$\lim_{t_1 \uparrow t_2} E|u(t_2,x)-u(t_1,x)| = 0, \text{ for a fixed } x \in \mathbf{R}^d.$$

We have the following a priori result about the continuity of the second moment of the solution to (1.1).

**Theorem 2.10** *Suppose that condition 2.2 holds, then for each $x \in \mathbf{R}$, the unique solution to (1.1) is mean square continuous in time. That is for each $x \in \mathbf{R}$, the function $t \to E[|u(t,x)|^2]$ is continuous.*

The compensated Poisson measure is a martingale-valued process and its second moment exists and hence the mean-square time continuity. We apply the Fourier transform via Plancherel's theorem to estimate the integrals.

Next, we present the time continuity of the solution to (1.2).

**Theorem 2.11** *Suppose that condition 2.5 holds, then for each $x \in \mathbf{R}^d$, the unique*

*solution to (1.2) is mean continuous in time. That is for each $x \in \mathbf{R}^d$, the function $t \to \mathrm{E}[|u(t,x)|]$ is continuous.*

For the non-compensated equation, it requires that we take the first moment of the mild solution. There are no useful tools to estimate the heat kernel on the integrals and we therefore make use of the heat kernel of alpha-stable processes; thus applying some explicit bounds on the fractional heat kernel to obtain more precise results.

## 3 METHODOLOGY

We present some preliminary concepts used in this paper. Let $(\Omega, \mathsf{F}, \mathsf{P})$ be a complete probability space and $(\mathbf{R}^d, \mathsf{B}(\mathbf{R}^d))$ be a measurable space. Let $\mathsf{M}$ be the space of all $\overline{\mathsf{Z}}_+ = \mathsf{Z}_+ \cup \{\infty\}$-valued measures on $(\mathbf{R}^d, \mathsf{B}(\mathbf{R}^d))$ and consider the measurable space $(\mathsf{M}, \mathsf{B}_\mathsf{M})$ with
$$\mathsf{B}_\mathsf{M} := \sigma(\nu(A) : A \in \mathsf{B}(\mathbf{R}^d)), \text{ for each } \nu \in \mathsf{M}.$$

**Definition 3.1** *(Poisson random measure) Let $\nu$ be a $\sigma$-finite measure on $(\mathbf{R}^d, \mathsf{B}(\mathbf{R}^d))$. A random variable $N : (\Omega, \mathsf{F}) \to (\mathsf{M}, \mathsf{B}_\mathsf{M})$ with intensity measure $\nu$ is called a Poisson random measure on $(\mathbf{R}^d, \mathsf{B}(\mathbf{R}^d))$ if the following conditions hold:*

• For all $A \in \mathsf{B}, N(A) : \Omega \to \overline{\mathsf{Z}}_+$ is Poisson distributed with parameter $\mathrm{E}[N(A)] = \nu(A)$, that is:
$$P[N(A) = n] = \frac{[\mathrm{E}[N(A)]]^n \exp(-\mathrm{E}[N(A)])}{n!}, n \in \mathsf{N} \cup \{0\}.$$
If $\mathrm{E}[N(A)] = +\infty$ then $N(A) = +\infty$ P-a.s.

• If $A_1, ..., A_k$ are pairwise disjoint then $N(A_1), ..., N(A_k)$ are independent.

Next we state some existence theorems for the Poisson random measure. The first one states that given a $\sigma$-finite measure on a space $X$, we can find or construct a Poisson random measure $N$ and it is given below.

**Theorem 3.2** *Given a $\sigma$-finite measure $\nu$ on $(X, \mathsf{B}(X))$, there exists a Poisson random measure $N$ such that $\mathrm{E}[N(A)] = \nu(A)$ for $A \in \mathsf{B}(X)$. If $N(A) = \infty$, then $\nu(A) = \infty$.*

*Proof.* The proof can be found in (Ikeda and Watanabe, 1981; Ikeda and Watanabe, 1989).

Let $(\Omega, \mathsf{F}, \mathsf{P})$ be a complete probability space and $(X, \mathsf{B}(X))$ a measurable space. Let the set of all point functions taking values in X, be denoted by $\Pi_X$ and $\mathsf{B}(\Pi_X)$ the smallest $\sigma$-algebra such that every mapping $p \to N_p((0,t] \times A)$ for all $A \in \mathsf{B}(X)$ is measurable, that's,
$$\mathsf{B}(\Pi_X) := \sigma(\Pi_X \ni p \mapsto N_p((0,t] \times A) : A \in \mathsf{B}(X), t > 0).$$

**Definition 3.3 (Point process)** *A Point process $p$ on $X$ is a $(\Pi_X, \mathsf{B}(\Pi_X))$-valued random variable. In other words, $p$ is defined on some probability space $(\Omega, \mathsf{F}, \mathsf{P})$, measurable and it spits out a point function from $\Pi_X$, that's, a random variable $p:(\Omega, \mathsf{F}) \to (\Pi_X, \mathsf{B}(\Pi_X))$.*

**Definition 3.4 (Poisson point process)** *A Point process $p$ is said to be a Poisson point process if the corresponding counting measure $N_p(\mathrm{d}t\mathrm{d}x)$ is a Poisson random measure on $(0, \infty) \times X$.*

The second existence theorem characterises a stationary Point process.

**Theorem 3.5** *Given a $\sigma$-finite measure $n(\mathrm{d}x)$ on $(X, \mathsf{B}(X))$, then there exists a stationary Poisson point process $p$ if the random measure $N_p(\mathrm{d}t, \mathrm{d}x)$ is of the form*
$$E[N_p(\mathrm{d}t, \mathrm{d}x)] = n_p(\mathrm{d}t, \mathrm{d}x) = \mathrm{d}tn(\mathrm{d}x).$$
*Proof.* (Ikeda and Watanabe, 1981).

Now applying the above theorem with $X := \mathbf{R}^d \times \mathbf{R}^d$ and $\mathsf{B}(X) := \mathsf{B}(\mathbf{R}^d) \otimes \mathsf{B}(\mathbf{R}^d)$. We will take $n(\mathrm{d}x, \mathrm{d}h) := \mathrm{d}x\nu(\mathrm{d}h)$. One set of the vectors will play the role of position while the other will play the role of "jumps". By the above theorem, we have a Poisson point process $p(s) \in \mathbf{R}^d \times \mathbf{R}^d$. The Poisson random measure is thus given by the following

$$N_p((0,t], A \times B) := \#\{s \leq t; s \in D_p; p(s) \in A \times B\}.$$

**Definition 3.6 (*Jump of a Lévy Process*)** *The jump process $\Delta X_t$ at time $t \geq 0$ is defined by $\Delta X_t := X_t - X_{t^-}$ where $X_{t^-}$ is the left limit of the process $X_t$ at the point $t$.*

**Definition 3.7 (*Jump measure*)** *Let $(\Delta X_t \neq 0, t > 0)$ be the jump process and the set $A \in \mathsf{B}(\mathbf{R}^d)$ bounded below, then one defines the jump measure by*
$$N(t, A) = \#\{0 \leq s \leq t : \Delta X_s \in A\} = \sum_{0 < s \leq t} I_A(\Delta X_s).$$
The jump measure counts the number of jumps of the process between $0$ and $t$ such that their sizes fall into $A$.

**Definition 3.8 (Compensated Poisson process)** *For a Poisson process, $N((0,t], A \times B)$ such that*
$$E[N((0,t], A \times B)] = t|A|\nu(B), \text{ for all } A, B \in \mathsf{B}(\mathbf{R}^d),$$
one defines the compensated Poisson process by
$$\tilde{N}((0,t], A \times B) := N((0,t], A \times B) - t|A|\nu(B),$$
for any $t > 0$ and any $A, B \in \mathsf{B}(\mathbf{R}^d)$ provided that $|A|\nu(B) < \infty$.

**Definition 3.9** *(Lévy measure) Let* $(\Omega, \mathsf{F}, \mathrm{P})$ *be a complete probability space. The measure $\nu$ defined by*

$$\nu(A) = \mathrm{E}[N((0,1], A)] = \mathrm{E}[\sum_{0 < s \leq 1} I_A(\Delta X_s)]$$

for all $A \in \mathsf{B}(\mathbf{R}^d)$ is said to be a Lévy measure of the process $X$ with E an expectation with respect to the measure $\mathrm{P}$. Suppose that $\nu$ is a Lévy measure on $\mathbf{R}^d$; then it satisfies the following

$$\int_{\mathbf{R}^d} (1 \wedge h^2) \nu(\mathrm{d}h) < \infty.$$

### 3.1 The Poisson Discontinuous Integrals

We now make sense of the discontinuous integrals. Let $(\Omega, \mathsf{F}, \{\mathsf{F}_t\}_{t \geq 0}, \mathrm{P})$ be a complete filtered Probability space and $(\mathbf{R}^d, \mathsf{B}(\mathbf{R}^d))$ be a measurable space. Let $\mathsf{F}_t$ be defined by

$$\mathsf{F}_t := \sigma(N_p([0,t], A \times B, \cdot) : A \times B \in \mathsf{B}(\mathbf{R}^d) \times \mathsf{B}(\mathbf{R}^d)) \vee \mathsf{N},$$

where $t > 0$ and $\mathsf{N}$ denotes the null set of $\mathsf{F}$. We can write the Poisson random measure as

$$N_p((0,t], A \times B) := \sum_{s \in D_p, s \leq t} I_{A \times B}(p_x(s), p_h(s)),$$

where we define $p(s) := (p_x(s), p_h(s))$.

Recall that in our case, we have $\mathrm{E}[N_p((0,t], A \times B)] = t|A|\nu(B)$. We now describe the stochastic integral with respect to this Poisson random measure. We will need to define the class of integrand precisely.

**Definition 3.10** (The non-compensated Integral)

$$H_p^1 := \{f(t,x,h) : f \text{ is } \{\mathsf{F}_t\} - predictable \text{ and } \int_0^t \int_{\mathbf{R}^d} \int_{\mathbf{R}^d} \mathrm{E} |f(s,x,h)| \mathrm{d}s \mathrm{d}x \nu(\mathrm{d}h) < \infty\}.$$

The following integral can now be defined for all $f \in H_p^1$:

$$\int_0^t \int_{\mathbf{R}^d} \int_{\mathbf{R}^d} f(s,x,h,.) N_p(\mathrm{d}s, \mathrm{d}x, \mathrm{d}h) = \sum_{s \leq t, s \in D_p} f(s, p_x(s), p_h(s))$$

as the a.s sum of the following absolutely convergent sum.

**Definition 3.11** (The compensated Integral) *Define, similarly, for $f$ satisfying the square-integrability condition*

$$H_p^2 = \{f(t,x,h) : f \text{ is } \{\mathsf{F}_t\} - predictable \text{ and } \int_0^t \int_{\mathbf{R}^d} \int_{\mathbf{R}^d} \mathrm{E} |f(s,x,h)|^2 \mathrm{d}s \mathrm{d}x \nu(\mathrm{d}h) < \infty\}.$$

Then for all $f \in H_p^2$, one defines the integral as follows

$$\int_0^t \int_{\mathbf{R}^d} \int_{\mathbf{R}^d} f(s,x,h)\tilde{N}_p(ds,dx,dh) = \sum_{s\leq t, s\in D_p} f(s,p_x(s),p_h(s))$$

$$- \int_0^t \int_{\mathbf{R}^d} \int_{\mathbf{R}^d} f(s,x,h)ds dx \nu(dh)$$

as the a.s sum of the following absolutely convergent sum.

### 3.2 Symmetric $\alpha$-stable processes

**Definition 3.12** *(Stable process) A random variable $X$ is said to be stable if there exist real valued sequences $(c_n, n \in \mathbb{N})$ and $(d_n, n \in \mathbb{N})$ with each $c_n > 0$ such that*

$$X_1 + X_2 + \ldots + X_n \stackrel{d}{=} c_n X + d_n \tag{3.1}$$

where $X_1 + X_2 + \ldots + X_n$ are independent copies of $X$. The random variable $X$ is said to be strictly stable if each $d_n = 0$. It has been shown (Feller, 1971) that the only choice of $c_n$ in (3.1) is of the form

$$c_n = \sigma n^{\frac{1}{\alpha}}, \quad 0 < \alpha \leq 2.$$

The parameter $\alpha$ plays a key role in the investigation of stable random variables and it is called the "index of stability". The operator $-(-\Delta)^{\alpha/2}$ is the fractional Laplacian of the $L^2$-generator of a symmetric stable process $X_t$ of order $\alpha$.

**Definition 3.13** *(Symmetric stable process) A symmetric $\alpha$-stable process $X$ on $\mathbf{R}^d$ is a Lévy process whose transition density $p(t,x)$ relative to Lebesgue measure is uniquely determined by its Fourier transform:*

$$\mathrm{E}[\exp(i\xi X_t)] = \int_{\mathbf{R}^d} e^{i\langle x,\xi \rangle} p(t,x) dx = e^{-t|\xi|^\alpha}, \quad \xi \in \mathbf{R}^d.$$

We present some required properties of $p(t,x)$ which come in handy in the proof of our results, see (Sugitani, 1975).

$$p(t,x) = t^{-d/\alpha} p(1, t^{-1/\alpha} x)$$
$$p(st,x) = t^{-d/\alpha} p(s, t^{-1/\alpha} x). \tag{3.2}$$

From the above relation, $p(t,0) = t^{-d/\alpha} p(1,0)$, is a decreasing function of $t$. The heat kernel $p(t,x)$ is also a decreasing function of $|x|$, that's

$$|x| \geq |y| \text{ implies that } p(t,x) \leq p(t,y).$$

This and equation (3.2) imply that for all $t \geq s$,

$$p(t,x) = p(t,|x|) = p(s \cdot \frac{t}{s}, |x|) = (\frac{t}{s})^{-d/\alpha} p(s, (\frac{t}{s})^{-1/\alpha} |x|)$$

$$\geq (\frac{s}{t})^{d/\alpha} p(s, |x|) \quad (since \ (\frac{t}{s})^{-1/\alpha} |x| \leq |x|)$$

and $\quad p(t,x) \geq (\frac{s}{t})^{d/\alpha} p(s,x).$

**Proposition 3.14** *Let $p(t,x)$ be the transition density of a strictly $\alpha$-stable process. If $p(t,0) \leq 1$ and $a \geq 2$, then*

$$p(t, \frac{1}{a}(x-y)) \geq p(t,x)p(t,y) \, \forall x, y \in \mathbf{R}^d.$$

*Proof.* Given that

$$\frac{1}{a}|x-y| \leq \frac{2}{a}|x| \vee \frac{2}{a}|y| \leq |x| \vee |y|,$$

then it follows from the above that,

$$p(t, \frac{1}{a}(x-y)) \geq p(t, |x| \vee |y|) \geq p(t, |x|) \wedge p(t, |y|) \geq p(t, |x|)p(t, |y|) = p(t,x)p(t,y).$$

The transition density also satisfies the following Chapman-Kolmogorov equation,

$$\int_{\mathbf{R}^d} p(t,x)p(s,x)\mathrm{d}x = p(t+s,0).$$

**Lemma 3.15** *Suppose that $p(t,x)$ denotes the heat kernel for a strictly stable process of order $\alpha$. Then the following estimate holds.*

$$p(t,x,y) \approx t^{-d/\alpha} \wedge \frac{t}{|x-y|^{d+\alpha}} \quad \text{for all} \quad t > 0 \quad \text{and} \quad x, y \in \mathbf{R}^d.$$

Here and in the sequel, for two non-negative functions $f, g$, $f \approx g$ means that there exists a positive constant $c > 1$ such that $c^{-1}g \leq f \leq cg$ on their common domain of definition.

## 4 MAIN RESULTS

We present the proofs to our main results here.

### 4.1 Proof of mean-square continuity in time of theorem 2.10

The mild solution is given by

$$u(t,x) = \int_{\mathbf{R}} p(t,x,y)u(0,y)\mathrm{d}y + \lambda \int_0^t \int_{\mathbf{R}} \int_{\mathbf{R}} p(t-s,x,y)\sigma(u(s,y),h)\tilde{N}(\mathrm{d}h,\mathrm{d}y,\mathrm{d}s).$$

We assume $0 < t_1 < t_2$, then for fixed $x \in \mathbf{R}$

$$u(t_2,x) - u(t_1,x) = \int_{\mathbf{R}} [p(t_2,x,y) - p(t_1,x,y)]u(0,y)\mathrm{d}y$$

$$+ \lambda \int_{\mathbf{R}} \int_{\mathbf{R}} \int_0^{t_1} [p(t_2-s,x,y) - p(t_1-s,x,y)]\sigma(u(s,y),h)\tilde{N}(\mathrm{d}h,\mathrm{d}y,\mathrm{d}s)$$

$$+ \lambda \int_{\mathbf{R}} \int_{\mathbf{R}} \int_{t_1}^{t_2} p(t_2-s,x,y)\sigma(u(s,y),h)\tilde{N}(\mathrm{d}h,\mathrm{d}y,\mathrm{d}s).$$

We make the following definitions,

$$D_0 = \int_{\mathbf{R}} [p(t_2,x,y) - p(t_1,x,y)]u(0,y)\mathrm{d}y$$

$$D_1 = \lambda \int_{\mathbf{R}} \int_{\mathbf{R}} \int_0^{t_1} [p(t_2-s,x,y) - p(t_1-s,x,y)]\sigma(u(s,y),h)\tilde{N}(\mathrm{d}h,\mathrm{d}y,\mathrm{d}s)$$

$$D_2 = \lambda \int_{\mathbf{R}} \int_{\mathbf{R}} \int_{t_1}^{t_2} p(t_2 - s, x, y) \sigma(u(s, y), h) \tilde{N}(dh, dy, ds).$$

The proof of theorem 2.10 will be a consequence of the following Lemma(s).

**Lemma 4.1** *For all $\beta > 0$, $0 < t_1 < t_2$ and $x \in \mathbf{R}$,*

$$|D_0|^2 \leq \frac{c_0}{2\pi} \int_{\mathbf{R}} e^{-2t_1 \mathrm{Re}\Psi(\xi)} |1 - e^{-(t_2-t_1)\Psi(\xi)}|^2 \, d\xi.$$

*Proof.* We start by writing

$$\mathrm{E}|D_0|^2 = |D_0|^2 = \left| \int_{\mathbf{R}} [p(t_2, x, y) - p(t_1, x, y)] u(0, y) dy \right|^2$$

$$\leq \int_{\mathbf{R}} |u(0, y)|^2 \, dy \int_{\mathbf{R}} |p(t_2, x, y) - p(t_1, x, y)|^2 \, dy$$

$$\leq c_0 \int_{\mathbf{R}} |p(t_2, x, y) - p(t_1, x, y)|^2 \, dy = c_0 \|p(t_2, .) - p(t_1, .)\|_{L^2(\mathbf{R})}^2.$$

By Plancherel's theorem:

$$\|p(t_2, .) - p(t_1, .)\|_{L^2(\mathbf{R})}^2 = \|\hat{p}(t_2, .) - \hat{p}(t_1, .)\|_{L^2(\mathbf{R})}^2 = \frac{1}{2\pi} \int_{\mathbf{R}} e^{-2t_1 \mathrm{Re}\Psi(\xi)} |1 - e^{-(t_2-t_1)\Psi(\xi)}|^2 \, d\xi.$$

Therefore,

$$\mathrm{E}|D_0|^2 \leq \frac{c_0}{2\pi} \int_{\mathbf{R}} e^{-2t_1 \mathrm{Re}\Psi(\xi)} |1 - e^{-(t_2-t_1)\Psi(\xi)}|^2 \, d\xi.$$

**Lemma 4.2** *For all $\beta > 0$, $0 < t_1 < t_2$ and $x \in \mathbf{R}$,*

$$\mathrm{E}|D_1|^2 \leq \frac{\lambda^2 \mathrm{KLip}_\sigma^2}{2\pi} \|u\|_{2,\beta}^2 e^{\beta t_1} \int_{\mathbf{R}} \frac{|1 - e^{-(t_2-t_1)\Psi(\xi)}|^2}{\beta + 2\mathrm{Re}\Psi(\xi)} \, d\xi.$$

*Proof.* By Itô's isometry, we obtain

$$\mathrm{E}|D_1|^2 = \lambda^2 \int_{\mathbf{R}} \int_{\mathbf{R}} \int_0^{t_1} |p(t_2 - s, x, y) - p(t_1 - s, x, y)|^2 \, \mathrm{E}|\sigma(u(s, y), h)|^2 \, \nu(dh) dy ds$$

$$\leq \lambda^2 \mathrm{KLip}_\sigma^2 \int_{\mathbf{R}} \int_0^{t_1} |p(t_2 - s, x, y) - p(t_1 - s, x, y)|^2 \, \mathrm{E}|u(s, y)|^2 \, dy ds$$

$$\leq \lambda^2 \mathrm{KLip}_\sigma^2 \|u\|_{2,\beta}^2 \int_0^{t_1} e^{\beta s} \|\hat{p}(t_2 - s, .) - \hat{p}(t_1 - s, .)\|_{L^2(\mathbf{R})}^2 \, ds.$$

But

$$\|\hat{p}(t_2 - s, .) - \hat{p}(t_1 - s, .)\|_{L^2(\mathbf{R})}^2 = \frac{1}{2\pi} \int_{\mathbf{R}} e^{-2(t_1-s)\mathrm{Re}\Psi(\xi)} |1 - e^{-(t_2-t_1)\Psi(\xi)}|^2 \, d\xi.$$

Therefore,

$$\mathrm{E}|D_1|^2 \leq \frac{\lambda^2 \mathrm{KLip}_\sigma^2}{2\pi} \|u\|_{2,\beta}^2 \int_{\mathbf{R}} d\xi \frac{|1 - e^{-(t_2-t_1)\Psi(\xi)}|^2}{\beta + 2\mathrm{Re}\Psi(\xi)} [1 - e^{-t_1(\beta + 2\mathrm{Re}\Psi(\xi))}]$$

$$\leq \frac{\lambda^2 \mathrm{KLip}_\sigma^2}{2\pi} \|u\|_{2,\beta}^2 e^{\beta t_1} \int_{\mathbf{R}} \frac{|1 - e^{-(t_2-t_1)\Psi(\xi)}|^2}{\beta + 2\mathrm{Re}\Psi(\xi)} \, d\xi.$$

**Lemma 4.3** *For all $\beta > 0$, $0 < t_1 < t_2$ and $x \in \mathbf{R}$,*

$$E|D_2|^2 \leq \frac{\lambda^2 \mathrm{KLip}_\sigma^2}{2\pi} \|u\|_{2,\beta}^2 \int_\mathbf{R} \frac{d\xi}{\beta + 2\mathrm{Re}\Psi(\xi)} \cdot e^{\beta t_2} [1 - e^{-(t_2 - t_1)(\beta + 2\mathrm{Re}\Psi(\xi))}].$$

*Proof.* Take second moment of the solution

$$E|D_2|^2 = \lambda^2 \int_\mathbf{R} \int_\mathbf{R} \int_{t_1}^{t_2} |p(t_2 - s, x, y)|^2 \, E|\sigma(u(s,y), h)|^2 \, \nu(dh) dy ds$$

$$\leq \lambda^2 \mathrm{KLip}_\sigma^2 \|u\|_{2,\beta}^2 \int_{t_1}^{t_2} e^{\beta s} \|\hat{p}(t_2 - s, .)\|_{L^2(\mathbf{R})}^2 \, ds$$

$$\leq \frac{\lambda^2 \mathrm{KLip}_\sigma^2}{2\pi} \|u\|_{2,\beta}^2 \int_{t_1}^{t_2} ds \, e^{\beta s} \int_\mathbf{R} e^{-2(t_2 - s)\mathrm{Re}\Psi(\xi)} d\xi$$

$$= \frac{\lambda^2 \mathrm{KLip}_\sigma^2}{2\pi} \|u\|_{2,\beta}^2 \int_\mathbf{R} \frac{d\xi}{\beta + 2\mathrm{Re}\Psi(\xi)} \cdot e^{\beta t_2} [1 - e^{-(t_2 - t_1)(\beta + 2\mathrm{Re}\Psi(\xi))}].$$

*Proof of Theorem 2.10.* Combining Lemma 4.1, 4.2 and 4.3, therefore

$$E|u(t_2, x) - u(t_1, x)|^2 \leq \frac{C}{2\pi} \int_\mathbf{R} e^{-2t_1 \mathrm{Re}\Psi(\xi)} |1 - e^{-(t_2 - t_1)\Psi(\xi)}|^2 \, d\xi$$

$$+ \frac{\lambda^2 \mathrm{KLip}_\sigma^2}{2\pi} \|u\|_{2,\beta}^2 \, e^{\beta t_1} \int_\mathbf{R} \frac{|1 - e^{-(t_2 - t_1)\Psi(\xi)}|^2}{\beta + 2\mathrm{Re}\Psi(\xi)} \, d\xi$$

$$+ \frac{\lambda^2 \mathrm{KLip}_\sigma^2}{2\pi} \|u\|_{2,\beta}^2 \int_\mathbf{R} \frac{d\xi}{\beta + 2\mathrm{Re}\Psi(\xi)} \, e^{\beta t_2} [1 - e^{-(t_2 - t_1)(\beta + 2\mathrm{Re}\Psi(\xi))}].$$

Then

$$\lim_{\delta \downarrow 0} \sup_{|t_1 - t_2| < \delta} E|u(t_2, x) - u(t_1, x)|^2 \leq 0$$

and therefore

$$\lim_{t_1 \uparrow t_2} E|u(t_2, x) - u(t_1, x)|^2 = 0 \text{ for a fixed } x \in \mathbf{R}.$$

### 4.2 FURTHER RESULTS

For the mean time continuity of solution to the non-compensated heat equation, we consider the case for $\alpha$-stable processes. This is because we have explicit estimates on the heat kernel for $\alpha$ stable processes. For the condition on the existence and uniqueness result for the stable process, we have:

**Theorem 4.4** *Suppose that $C_{d,\alpha,\beta} < \frac{1}{\lambda \mathrm{KLip}_\sigma}$ for positive constants K, $\mathrm{Lip}_\sigma$, then there exists a solution u that is unique up to modification.*

The proof of the above theorem is based on the following Lemma 4.5 and Lemma 4.6.

Now let

$$A^\alpha u(t,x) := \lambda \int_0^t \int_{\mathbf{R}^d} \int_{\mathbf{R}^d} p^\alpha(t-s,x,y)\sigma(u(s,y),h)N(dh,dy,ds),$$

and then the following Lemma(s):

**Lemma 4.5** *Suppose that $u$ is predictable and $\|u\|_{1,\beta} < \infty$ for all $\beta > 0$ and $\sigma(u,h)$ satisfies condition 2.5, then*

$$\|A^\alpha u\|_{1,\beta} \le C_{d,\alpha,\beta}\lambda K[1+\text{Lip}_\sigma \|u\|_{1,\beta}],$$

where $C_{d,\alpha,\beta} := \dfrac{2C(d,\alpha)d+\alpha}{d+\alpha-1}\dfrac{\Gamma(\gamma+1)}{\beta^{\gamma+1}}.$

*Proof.* Taking first moment of the solution, we have

$$E|A^\alpha u(t,x)| = \lambda \int_0^t \int_{\mathbf{R}^d} \int_{\mathbf{R}^d} |p^\alpha(t-s,x,y)| E|\sigma(u(s,y),h)| \nu(dh)dyds$$

$$\le \lambda K \int_0^t \int_{\mathbf{R}^d} |p^\alpha(t-s,x,y)|[1+\text{Lip}_\sigma E|u(s,y)|]dyds.$$

Next, Multiply through by $\exp(-\beta t)$, to get

$$e^{-\beta t}E|A^\alpha u(t,x)| \le \lambda K \int_0^t \int_{\mathbf{R}^d} e^{-\beta(t-s)} |p^\alpha(t-s,x,y)|\{e^{-\beta s}[1+\text{Lip}_\sigma E|u(s,y)|]\}dyds$$

$$\le \lambda K \text{Lip}_\sigma \sup_{s\ge 0}\sup_{y\in\mathbf{R}^d}\{e^{-\beta s}[1+\text{Lip}_\sigma E|u(s,y)|]\}\int_0^t \int_{\mathbf{R}^d} e^{-\beta(t-s)}|p^\alpha(t-s,x,y)|dyds.$$

Then we obtain that

$$\|A^\alpha u\|_{1,\beta} \le \lambda K[1+\text{Lip}_\sigma \|u\|_{1,\beta}]\sup_{t\ge 0}\int_0^t \int_{\mathbf{R}^d} e^{-\beta(t-s)} |p^\alpha(t-s,x,y)|dyds$$

$$\le \lambda K[1+\text{Lip}_\sigma \|u\|_{1,\beta}]\int_0^\infty \int_{\mathbf{R}^d} e^{-\beta s}|p^\alpha(s,y)|dyds$$

$$\le \lambda K[1+\text{Lip}_\sigma \|u\|_{1,\beta}]\int_0^\infty \int_{\mathbf{R}^d} e^{-\beta s}\{C(\frac{s}{|y|^{d+\alpha}}\wedge s^{-\frac{d}{\alpha}})\}dyds.$$

The last inequality follows by Lemma 3.15. Let's assume that $\dfrac{s}{|y|^{d+\alpha}} \le s^{-\frac{d}{\alpha}}$ which holds only when $|y|^\alpha \ge s$. Therefore

$$\|A^\alpha u\|_{1,\beta} \le C(d,\alpha)\lambda K[1+\text{Lip}_\sigma \|u\|_{1,\beta}]\int_0^\infty dse^{-\beta s}[s\int_{|y|\ge s^{1/\alpha}}\frac{dy}{|y|^{d+\alpha}}$$

$$+ s^{-d/\alpha}\int_{|y|<s^{1/\alpha}}dy]$$

$$= C(d,\alpha)\lambda K[1+\text{Lip}_\sigma \|u\|_{1,\beta}]\int_0^\infty dse^{-\beta s}[s(-\int_{-\infty}^{s^{1/\alpha}} y^{-(d+\alpha)}dy$$

$$+ \int_{s^{1/\alpha}}^\infty y^{-(d+\alpha)}dy) + 2s^{(1-d)/\alpha}]$$

$$= C(d,\alpha)\lambda K[1+\text{Lip}_\sigma \|u\|_{1,\beta}]\int_0^\infty dse^{-\beta s}[s(-\frac{y^{-(d+\alpha-1)}}{1-d-\alpha}\Big|_{-\infty}^{s^{1/\alpha}}$$

$$+ \frac{y^{-(d+\alpha-1)}}{1-d-\alpha}\Big|_{s^{1/\alpha}}^{\infty}) + 2s^{(1-d)/\alpha}]$$

$$= C(d,\alpha)\lambda K[1+\text{Lip}_\sigma \|u\|_{1,\beta}]\int_0^\infty ds e^{-\beta s}[s(-\frac{2}{1-d-\alpha}s^{(1-d-\alpha)/\alpha})$$

$$+ 2s^{(1-d)/\alpha}]$$

$$= C(d,\alpha)\lambda K[1+\text{Lip}_\sigma \|u\|_{1,\beta}]\int_0^\infty ds e^{-\beta s}[\frac{2}{d+\alpha-1}s^{1+(1-d-\alpha)/\alpha}$$

$$+ 2s^{(1-d)/\alpha}].$$

Thus

$$\|A^\alpha u\|_{1,\beta} \leq 2C(d,\alpha)\lambda K[1+\text{Lip}_\sigma \|u\|_{1,\beta}]\frac{d+\alpha}{d+\alpha-1}\int_0^\infty s^\gamma e^{-\beta s}ds,$$

where $\gamma = (1-d)/\alpha$. Then we compute the integral $I_{\beta,\gamma} := \int_0^\infty s^\gamma e^{-\beta s}ds$. Let $\tau = \beta s, ds = \frac{d\tau}{\beta}$, now therefore,

$$I_{\beta,\gamma} = \frac{1}{\beta^{\gamma+1}}\int_0^\infty \tau^\gamma e^{-\tau}d\tau = \frac{1}{\beta^{\gamma+1}}\int_0^\infty \tau^{(\gamma+1)-1} e^{-\tau}d\tau = \frac{\Gamma(\gamma+1)}{\beta^{\gamma+1}}.$$

Therefore

$$\|A^\alpha u\|_{1,\beta} \leq 2C(d,\alpha)\lambda K[1+\text{Lip}_\sigma \|u\|_{1,\beta}][\frac{d+\alpha}{d+\alpha-1}\frac{\Gamma(\gamma+1)}{\beta^{\gamma+1}}].$$

**Lemma 4.6** *Suppose $u$ and $v$ are two predictable random field solutions satisfying $\|u\|_{1,\beta} + \|v\|_{1,\beta} < \infty$ for all $\beta > 0$ and $\sigma(u,h)$ satisfies condition 2.5, then*

$$\|A^\alpha u - A^\alpha v\|_{1,\beta} \leq C_{d,\alpha,\beta}\lambda K \text{Lip}_\sigma \|u-v\|_{1,\beta}.$$

*Proof.* Similar steps as Lemma 4.2

### 4.3 Proof of mean continuity in time of theorem 2.2

The solution is given by

$$u(t,x) = \int_{\mathbf{R}^d} p^\alpha(t,x,y)u(0,y)dy$$

$$+ \lambda \int_0^t \int_{\mathbf{R}^d}\int_{\mathbf{R}^d} p^\alpha(t-s,x,y)\sigma(u(s,y),h)N(dh,dy,ds).$$

We assume $0 < t_1 < t_2$, then for fixed $x \in \mathbf{R}^d$

$$u(t_2,x) - u(t_1,x) = \int_{\mathbf{R}^d}[p^\alpha(t_2,x,y) - p^\alpha(t_1,x,y)]u(0,y)dy$$

$$+ \lambda \int_{\mathbf{R}^d}\int_{\mathbf{R}^d}\int_0^{t_1}[p^\alpha(t_2-s,x,y) - p^\alpha(t_1-s,x,y)]$$

$$\times \sigma(u(s,y,.),h)N(dh,dy,ds)$$

$$+ \lambda \int_{\mathbf{R}^d}\int_{\mathbf{R}^d}\int_{t_1}^{t_2} p^\alpha(t_2-s,x,y)\sigma(u(s,y),h)N(dh,dy,ds).$$

We make the following definitions,

$$D_3 = \int_{\mathbf{R}^d} [p^\alpha(t_2, x, y) - p^\alpha(t_1, x, y)] u(0, y) dy$$

$$D_4 = \lambda \int_{\mathbf{R}^d} \int_{\mathbf{R}^d} \int_0^{t_1} [p^\alpha(t_2 - s, x, y) - p^\alpha(t_1 - s, x, y)] \sigma(u(s, y), h) N(dh, dy, ds)$$

$$D_5 = \lambda \int_{\mathbf{R}^d} \int_{\mathbf{R}^d} \int_{t_1}^{t_2} p^\alpha(t_2 - s, x, y) \sigma(u(s, y), h) N(dh, dy, ds).$$

The proof of theorem 2.11 is a consequence of the following lemmas.

**Lemma 4.7** *For all $\beta > 0$, $0 < t_1 < t_2$, $x \in \mathbf{R}^d$ then*

$$|D_3| \leq 2c_0 C(d, \alpha) \frac{d + \alpha}{d + \alpha - 1} (t_2^{(1-d)/\alpha} - t_1^{(1-d)/\alpha})$$

*Proof.* Write,

$$\mathrm{E}|D_3| = |D_3| = |\int_{\mathbf{R}^d} [p^\alpha(t_2, x, y) - p^\alpha(t_1, x, y)] u(0, y) dy|$$

$$\leq \sup_{y \in \mathbf{R}^d} |u(0, y)| \int_{\mathbf{R}^d} |p^\alpha(t_2, x, y) - p^\alpha(t_1, x, y)| dy$$

$$= c_0 \int_{\mathbf{R}^d} |p^\alpha(t_2, x, y) - p^\alpha(t_1, x, y)| dy.$$

Using the estimates on the heat kernel of stable processes:

$$p^\alpha(t_2, x - y) - p^\alpha(t_1, x - y) \equiv (t_2^{-d/\alpha} \wedge \frac{t_2}{|x - y|^{d+\alpha}}) - (t_1^{-d/\alpha} \wedge \frac{t_1}{|x - y|^{d+\alpha}})$$

Therefore,

$$\mathrm{E}|D_3| \leq c_0 \{\int_{\mathbf{R}^d} (t_2^{-d/\alpha} \wedge \frac{t_2}{|x - y|^{d+\alpha}}) dy - \int_{\mathbf{R}^d} (t_1^{-d/\alpha} \wedge \frac{t_1}{|x - y|^{d+\alpha}}) dy\}.$$

But

$$\int_{\mathbf{R}^d} (t_2^{-d/\alpha} \wedge \frac{t_2}{|x - y|^{d+\alpha}}) dy = t_2 \int_{|x-y| \geq t_2^{1/\alpha}} \frac{dy}{|x - y|^{d+\alpha}} + t_2^{-d/\alpha} \int_{|x-y| < t_2^{1/\alpha}} dy$$

$$= \frac{2C(d, \alpha)}{d + \alpha - 1} t_2^{1+(1-d-\alpha)/\alpha} + 2 t_2^{(1-d)/\alpha} = 2 C(d, \alpha) \frac{d + \alpha}{d + \alpha - 1} t_2^{(1-d)/\alpha}.$$

Doing same for the other integral on $t_1$, therefore

$$|D_3| \leq 2c_0 C(d + \alpha) \frac{d + \alpha}{d + \alpha - 1} (t_2^{(1-d)/\alpha} - t_1^{(1-d)/\alpha}).$$

**Lemma 4.8** *For all $\beta > 0$, $0 < t_1 < t_2$ and $x \in \mathbf{R}^d$,*

$$\mathrm{E}|D_4| \leq 2\lambda \mathrm{KLip}_\sigma C(d, \alpha) \|u\|_{1,\beta} \frac{d + \alpha}{d + \alpha - 1} (-e^{\beta t_2} \int_{t_2 - t_1}^{t_2} z^{(1-d)/\alpha} e^{-\beta z} dz$$

$$+ e^{\beta t_1} \int_0^{t_1} z^{(1-d)/\alpha} e^{-\beta z} dz).$$

*Proof.* We begin by writing

$$E\mid D_4\mid = \int_{\mathbf{R}^d}\int_{\mathbf{R}^d}\int_0^{t_1} \mid p^\alpha(t_2-s,x,y)-p^\alpha(t_1-s,x,y)\mid E\mid \sigma(u(s,y),h)\mid \nu(\mathrm{d}h)\mathrm{d}y\mathrm{d}s$$

$$\le \lambda \mathrm{KLip}_\sigma \int_{\mathbf{R}^d}\int_0^{t_1}\mid p^\alpha(t_2-s,x,y)-p^\alpha(t_1-s,x,y)\mid E\mid u(s,y)\mid \mathrm{d}y\mathrm{d}s$$

$$\le \lambda \mathrm{KLip}_\sigma \|u\|_{1,\beta}\int_{\mathbf{R}^d}\int_0^{t_1} e^{\beta s}\mid p^\alpha(t_2-s,x,y)-p^\alpha(t_1-s,x,y)\mid \mathrm{d}s\mathrm{d}y$$

$$\le \lambda \mathrm{KLip}_\sigma C(d,\alpha)\|u\|_{1,\beta}\int_0^{t_1}\mathrm{d}s\, e^{\beta s}\{\int_{\mathbf{R}^d}((t_2-s)^{-d/\alpha}\wedge \frac{t_2-s}{\mid x-y\mid^{d+\alpha}})\mathrm{d}y$$

$$-\int_{\mathbf{R}^d}((t_1-s)^{-d/\alpha}\wedge \frac{t_1-s}{\mid x-y\mid^{d+\alpha}})\mathrm{d}y\}.$$

Similarly as above, therefore,

$$E\mid D_4\mid \le 2\lambda \mathrm{KLip}_\sigma C(d,\alpha)\|u\|_{1,\beta}\int_0^{t_1}\mathrm{d}s\, e^{\beta s}\{\frac{d+\alpha}{d+\alpha-1}((t_2-s)^{(1-d)/\alpha}-(t_1-s)^{(1-d)/\alpha})\}$$

and the result follows.

**Lemma 4.9** *For all* $\beta>0$, $0<t_1<t_2$ *and* $x\in \mathbf{R}^d$,

$$E\mid D_5\mid \le 2\lambda \mathrm{KLip}_\sigma C(d,\alpha)\|u\|_{1,\beta}e^{\beta t_2}\frac{d+\alpha}{d+\alpha-1}\int_0^{t_2-t_1}z^{(1-d)/\alpha}e^{-\beta z}\mathrm{d}z.$$

*Proof.* Taking an expectation of the solution,

$$E\mid D_5\mid = \lambda \int_{\mathbf{R}^d}\int_{\mathbf{R}^d}\int_{t_1}^{t_2}\mid p^\alpha(t_2-s,x,y)\mid E\mid \sigma(u(s,y),h)\mid \nu(\mathrm{d}h)\mathrm{d}y\mathrm{d}s$$

$$\le \lambda \mathrm{KLip}_\sigma \|u\|_{1,\beta}\int_{\mathbf{R}^d}\int_{t_1}^{t_2}e^{\beta s}\mid p^\alpha(t_2-s,x,y)\mid \mathrm{d}y\mathrm{d}s$$

$$= \lambda \mathrm{KLip}_\sigma C(d,\alpha)\|u\|_{1,\beta}\int_{t_1}^{t_2}\mathrm{d}s\, e^{\beta s}[(t_2-s)\int_{\mid x-y\mid \ge (t_2-s)^{1/\alpha}}\frac{\mathrm{d}y}{\mid x-y\mid^{d+\alpha}}$$

$$+(t_2-s)^{-d/\alpha}\int_{\mid x-y\mid <(t_2-s)^{1/\alpha}}\mathrm{d}y]$$

$$= \lambda \mathrm{KLip}_\sigma C(d,\alpha)\|u\|_{1,\beta}[\frac{2}{d+\alpha-1}\int_{t_1}^{t_2}e^{\beta s}(t_2-s)^{1+(1-d-\alpha)/\alpha}\mathrm{d}s$$

$$+2\int_{t_1}^{t_2}e^{\beta s}(t_2-s)^{(1-d)/\alpha}\mathrm{d}s]$$

$$= 2\lambda \mathrm{KLip}_\sigma C(d,\alpha)\|u\|_{1,\beta}e^{\beta t_2}\frac{d+\alpha}{d+\alpha-1}\int_0^{t_2-t_1}z^{(1-d)/\alpha}e^{-\beta z}\mathrm{d}z.$$

*Proof of Theorem 2.11.* Combining Lemma 4.7, 4.8 and 4.9, therefore

$$E\mid u(t_2,x)-u(t_1,x)\mid \le 2c_0 C(d,\alpha)\frac{d+\alpha}{d+\alpha-1}(t_2^{(1-d)/\alpha}-t_1^{(1-d)/\alpha})$$

$$+2\lambda \mathrm{KLip}_\sigma C(d,\alpha)\|u\|_{1,\beta}\frac{d+\alpha}{d+\alpha-1}(-e^{\beta t_2}\int_{t_2-t_1}^{t_2}z^{(1-d)/\alpha}e^{-\beta z}\mathrm{d}z$$

$$+e^{\beta t_1}\int_0^{t_1}z^{(1-d)/\alpha}e^{-\beta z}\mathrm{d}z)$$

$$+ 2\lambda \text{KLip}_\sigma C(d,\alpha)\|u\|_{1,\beta}\, e^{\beta t_2}\, \frac{d+\alpha}{d+\alpha-1} \int_0^{t_2-t_1} z^{(1-d)/\alpha} e^{-\beta z}\, dz\,.$$

Then

$$\lim_{\delta \downarrow 0} \sup_{|t_1-t_2|<\delta} E\,|u(t_2,x,.) - u(t_1,x,.)| \leq 0$$

and therefore

$$\lim_{\delta \downarrow 0} \sup_{|t_1-t_2|<\delta} E\,|u(t_2,x) - u(t_1,x)| = 0 \ \ for\ a\ fixed\ \ x \in \mathbf{R}^d.$$

## 5 DISCUSSION

Continuity property of the solutions to the discontinuous integral equations is one property or behaviour that the system should exhibit to guarantee (or in order to discuss) its global existence (or non-existence) under some non-linear conditions on $\sigma$. Our results show that the weak-predictable-random-field solutions are mean continuous and mean-square continuous in time for any time interval $t_2 \geq t_1 \geq 0$.